\documentclass[11pt]{amsart}
\begin{document}

\numberwithin{equation}{section}

\def\R{{\mathbb R}}
\def\K{{\mathcal K}}
\def\L{{\mathcal L}}
\def\O{{\Omega}}
\def\I{{\mathcal I}}
\def\T{{\mathcal T}}
\def\vep{{\varepsilon}}
\def\p{{\partial}}
\def\a{{\alpha}}
\def\b{{\beta}}
\def\e{{\eta}}
\def\g{{\gamma}}
\def\s{{\sigma}}
\def\t{{\tau}}
\def\E{{\mathbb E}}
\def\N{{\mathbb N}}
\def\F{{\mathcal F}}
\def\Pr{{\mathbb P}}
\def\EE{{\mathcal E}}
\def\o{{\omega}}
\def\l{{\ell}}
\def\vp{{\varphi}}
\def\vpo{{\varphi_{c_0}}}
\def\ve{{v^{\vep}}}
\def\un{{\mathrm{1~\hspace{-1.4ex}l}}}
\def\ue{{u^{\vep}}}
\def\pe{{\varphi_{c^{\varepsilon}}}}
\def\pet{{\varphi_{c^{\varepsilon}(t)}}}
\def\xe{{x^{\varepsilon}}}
\def\xet{{x^{\varepsilon}(t)}}
\def\ee{{\eta^{\varepsilon}}}
\def\ce{{c^{\varepsilon}}}
\def\cet{{c^{\varepsilon}(t)}}
\def\to{{\rightarrow}}
\def\te{{\tau^{\varepsilon}}}
\def\tea{{\tau^{\varepsilon}_{\alpha}}}
\def\ye{{y^{\varepsilon}}}
\def\ze{{z^{\varepsilon}}}
\def\be{{b^{\varepsilon}}}
\def\ae{{a^{\varepsilon}}}
\def\peto{{\varphi_{c^{\varepsilon}(\tau)}}}
\def\xeto{{x^{\varepsilon}(\tau)}}
\def\ceto{{c^{\varepsilon}(\tau)}}
\def\eeto{{\eta^{\varepsilon}(\tau)}}
\def\dt#1{\Frac{\partial#1}{\partial t}}
\def\dx#1{\Frac{\partial#1}{\partial x}}
\def\d3x#1{\Frac{\partial^3#1}{\partial x^3}}
\def\H#1#2{\widetilde H^{#1,#2}}
\def\Xbs{X_{b,s}}
\def\Xb0{X_{b,0}}
\def\X12{\widetilde X_{b,s_1,s_2}}
\def\Xs3{\widetilde X_{b,s,-3/8}}
\def\B03{\widetilde X_{b,0,-3/8}}
\def\xs{(1+|\xi|)^{2s}}

\newcommand{\Frac}{\displaystyle \frac}
\newcommand{\Sum}{\displaystyle \sum}
\newcommand{\Int}{\displaystyle \int}
\newcommand{\Sup}{\displaystyle \sup}

\newtheorem{Theorem}{Theorem}[section]
\newtheorem{Definition}[Theorem]{Definition}
\newtheorem{Proposition}[Theorem]{Proposition}
\newtheorem{Lemma}[Theorem]{Lemma}
\newtheorem{Corollary}[Theorem]{Corollary}
\newtheorem{Remark}[Theorem]{Remark}
\newtheorem{Example}[Theorem]{Example}

\title[Soliton dynamics for multiplicative KdV equation]
{Soliton dynamics for the
Korteweg-de Vries equation with multiplicative homogeneous noise}

\author[A. de Bouard]{Anne de Bouard$^{\small 1}$}
\author[A. Debussche]{Arnaud Debussche$^{\small 2}$}



\maketitle

\begin{center} \small
$^1$ Centre de Math\'ematiques Appliqu\'ees\\
UMR 7641, CNRS et Ecole Polytechnique, \\
Route de Saclay, \\
91128 PALAISEAU CEDEX, FRANCE
\end{center}

\vskip 0.1 in

\begin{center} \small
$^2$ IRMAR, ENS Cachan Bretagne, CNRS, UEB  \\
Av. Robert Schuman, \\
F-35170 BRUZ, FRANCE 
\end{center}

\vskip 0.1 in

\noindent
{\bf Abstract}
We consider a randomly perturbed Korteweg-de Vries equation. The perturbation is a random potential
depending both on space and time, with a white noise behavior in time, and a regular, but stationary
behavior in space. We investigate the dynamics of the soliton of the KdV equation in the presence of this
random perturbation, assuming that the amplitude of the perturbation is small. We estimate precisely the
exit time of the perturbed solution from a neighborhood of the modulated soliton, and we obtain the modulation
equations for the soliton parameters. We moreover prove a central limit theorem for the dispersive part of the
solution, and investigate the asymptotic behavior in time of the limit process.

\vskip 0.1 in

\section{Introduction}

Our aim is to describe the dynamics of a soliton solution of
the Korteweg-de Vries equation in the presence of a random potential, depending both
on space and time and  which is white in time. After the first paper \cite{Wa} showing ``superdiffusion"
of the soliton of the KdV equation in the presence of an external force which is a white noise in time
(see also \cite{BKKS}, \cite{He}),
the interest in such questions of soliton dynamics in
the presence of either deterministic or random perturbations has recently increased in the
mathematical community. In \cite{Ga}, e.g. the question is investigated with the help of
inverse scattering methods, for different types of time-white noise perturbations, still for
the KdV equation, while in
\cite{jfgs0}, \cite{jfgs1}, the case of a soliton of the NLS equation is studied, with the
presence of a slowly varying deterministic external potential. Random potential perturbations for NLS equations
have also been considered in \cite{Ga0} and \cite{dBF}. The diffusion of solitons of the KdV equation in the
presence of additive noise was numerically investigated in \cite{Pr}.
Also, in \cite{dBD1}, we studied the soliton dynamics for a KdV equation with an additive space-time noise.
Our aim here is to reproduce the analysis of \cite{dBD1} in the case of a random potential, which is stationary
in space : the solution of the stochastic equation starting from a soliton at initial time will then stay
close to a modulated soliton up to times small compared to  $\vep^{-2}$ where $\vep$ is the amplitude of
the random perturbation (see below).
In the present case, where the noise is multiplicative (the random potential) we are then able to analyze more
precisely the modulation equations for the soliton parameters and the linearized equation for the remaining
(dispersive) part
of the solution, and especially its asymptotic behavior in time.

We consider a stochastic KdV equation which may be written in It\^o form as
\begin{equation}
\label{1}
du+(\p_x^3 u+\frac12 \p_x(u^2)) dt = \vep u dW
\end{equation}
where $\vep>0$ is a small parameter, $u$ is a random process defined on
$(t,x) \in \R^+\times \R$,
$W$ is a Wiener process on $L^2(\R)$ whose covariance operator $\phi \phi^*$
is such that $\phi$ is a convolution
operator on $L^2(\R)$ defined by
$$ \phi f(x)= \Int_{\R} k(x-y) f(y) dy, \;  \mathrm{for} \; f \in L^2(\R). $$
The convolution kernel $k$ satisfies
\begin{equation}
\label{h1}
\|k\|_1:= \int_{\R} (k^2+(k')^2) dx < +\infty.
\end{equation}

Considering a complete orthonormal system $(e_i)_{i \in \N}$ in $L^2(\R)$,
we may alternatively write $W$ as
\begin{equation}
\label{2}
W(t,x)= \Sum_{i \in \N} \beta_i(t) \phi e_i(x),
\end{equation}
$(\beta_i)_{i \in \N}$ being an independent family of real valued Brownian
motions.
The correlation  function of the process $ W$ is then given by
$$ \E(W(t,x) W(s,y))=c(x-y)(s\wedge t), \quad x,y  \in \R,
\quad s,t>0, $$
where
$$c(z)=\int_{\R} k(z+u)k(u) du. $$
The existence and uniqueness of solutions for stochastic KdV equations of the
type (\ref{1}) but with an additive noise have been studied in \cite{dBD0},
\cite{dBDT99}, \cite{dBDT04}. The multiplicative case with homogeneous noise as described
above was considered in \cite{dBD06}:
assuming, together with the above condition, that $k$ is an integrable function of $x \in \R$
allowed us to prove the global existence and
uniqueness  of solutions to equation (\ref{1}) in the energy space $H^1(\R)$,
that is in the space where both the mass
\begin{equation}
\label{masse}
m(u)=\Frac{1}{2} \int_{\R} u^2(x) dx
\end{equation}
and the energy
\begin{equation}
\label{hamiltonien}
H(u)=\Frac{1}{2}\Int_{\R} (\p_x u)^2 dx -\Frac{1}{6} \Int_{\R} u^3 dx
\end{equation}
are well defined. Note that $m$ and $H$ are conserved for the
equation without noise, that is
\begin{equation}
\label{KdV}
\p_t u + \p_x^3 u + \frac12 \p_x (u^2) =0.
\end{equation}
Under the above conditions on $k$, it was then proved in \cite{dBD06} that for any given initial data
$u_0 \in H^1(\R)$, there is a unique solution $u$ of (\ref{1}) with paths a.s. continuous for $t\in \R^+$
with values in $H^1(\R)$. 

Our aim in this article is to analyze the qualitative influence of a noise on a
soliton solution of the deterministic equation.
More precisely, we study the qualitative behavior of solutions of (\ref{1}) in the limit
$\vep$ tends to zero, assuming that the initial state of the solution is a
soliton of
equation (\ref{KdV}). We recall indeed that equation (\ref{KdV}) possesses a
two-parameter family of solitary waves (or soliton) solutions, propagating with a constant velocity $c>0$, with the expression
$u_{c,x_0}(t,x)=\vp_c(x-ct+x_0)$, $x_0 \in \R$, where
\begin{equation}
\label{soliton}
\vp_c(x)= \Frac{3c}{2\cosh^2(\sqrt{c}\frac{x}{2})}
\end{equation}
satisfies the equation
\begin{equation}
\label{eqsoliton}
\vp_c^{\prime \prime} -c\vp_c +\frac12 \vp_c^2 =0.
\end{equation}

We do not recall here the well-known results concerning the stability of the soliton solutions
$u_{c,x_0}$ in equation (\ref{KdV}), but we refer to \cite{Ben}, \cite{BSS}, \cite{MM2}
or \cite{PW} for a review of the stability questions using PDE methods, or to \cite{GGKM} and
\cite{Sch} for a review of the stability of the solitons with the help of the inverse scattering transform.

Let us consider as in \cite{dBD1} the solution $\ue(t,x)$ of
equation (\ref{1}) which is such that $\ue(0,x)=\vp_0(x)$ where
$c_0>0$ is fixed. Then, in Section 2, we show, as we did in \cite{dBD1} for the additive equation
that up to times $C\vep^{-2}$, where $C$ is a constant, we may write the solution
$\ue$ as
\begin{equation}
\label{modulation}
\ue(t,x)= \pet (x-\xet) +\vep \ee (x-\xet)
\end{equation}
where the modulation parameters $\cet$ and $\xet$ satisfy a system of stochastic differential
equations and the remaining term $\vep \ee$ is small in $H^1(\R)$.
We then prove in Section 3 that the process $\ee$
converges as $\vep$ goes to zero, in quadratic mean, to a centered Gaussian
process $\eta$ which satisfies an additively driven linear equation, with a
conservative deterministic part;
we also investigate the behavior of the process $\eta$ as $t$ goes to
infinity and prove that $\e$ is in some sense a Ornstein-Uhlenbeck process,
with a unique Gaussian invariant measure.
In addition, the parameters $\xet$ and $\cet$ may be developed up to order one in $\vep$ and we get
$$ \left\{
\begin{array}{l}
d\xe=c_0dt+\vep B_1dt +\vep dB_2+o(\vep)\\
d\ce=\vep dB_1 +o(\vep),
\end{array}
\right.$$
where $B_1$ and $B_2$ are correlated real valued Brownian motions; keeping only the order one terms in those modulation
parameters, we then obtain a diffusion result on the modulated soliton similar to the result obtained by Wadati in \cite{Wa},
but with a different time exponent (see Section 4).

In all what follows,
$(.,.)$ will denote the inner product in $L^2(\R)$,
$$
(u,v)=\int_{\R} u(x)v(x) dx
$$
and we denote by
$\T_{x_0}$ the translation operator
defined for $\vp \in C(\R)$ by $(\T_{x_0} \vp) (x)=\vp(x+x_0)$.
Note that since the process $W$ is stationary in space, for any
$x_0 \in \R$ the process $\T_{x_0} W$ is still a Wiener process with
covariance $\phi \phi^*$. Indeed by (\ref{2}),
$$\T_{x_0}  W(t,x)=\sum_{k\in \N} (\phi e_k) (x+x_0) \b_k(t)=
\sum_{k\in \N} (\phi \tilde e_k)(x) \b_k(t),$$
with $\tilde e_k(x)=\T_{x_0}e_k$.

\section{Modulation and estimate on the exit time}
In this section, we prove the following theorem.

\begin{Theorem}
\label{t2}
Assume that the kernel $k$ of the noise satisfies (\ref{h1}) together with
$k\in L^1(\R)$
and let $c_0$ be fixed. For $\vep>0$, let $\ue(t,x)$,
as defined above, be the solution of (\ref{1}) with $u(0,x)=\vp_{c_0}(x)$.
Then there exists $\a_0 >0$ such that, for each $\a$, $0<\a\leq\a_0$, there
is a stopping time $\tea>0$ a.s. and there are semi-martingale processes $\cet$
and $\xet$, defined a.s. for $t\leq \tea$, with values respectively in
$\R^{+*}$ and $\R$, so that if we set $\vep \ee(t)=\ue(t,.+\xet)-\pet$, then
a.s. for $t\leq \tea$,
$
\|\vep \ee (t)\|_1 \leq \a
$
and
$
|\cet-c_0| \leq \a.
$
In addition, for $\a_0$ sufficiently small, and any
$\a\leq \a_0$, there is a constant $C>0$,
depending only on $\a$ and $c_0$, such that for any $T>0$,
there is an $\vep_0>0$, with, for each $\vep<\vep_0$,
\begin{equation}
\label{2.6}
\Pr(\tea\leq T) \leq
\exp\left(-\frac{C(\alpha,c_0)}{\vep^2T\|k\|_{H^1}^2}\right).
\end{equation}
\end{Theorem}

It was noticed heuristically in \cite{dBD1}, and proved in \cite{dBG} that in the
additive case, the use of the modulation parameters $\xet$ and $\cet$ was
necessary in order to get the estimate (\ref{2.6}).
Indeed, it was proved in \cite{dBG} that if we denote by
${\tilde \tau}^{\varepsilon,n}_{\alpha} =\inf \{ t>0, \|u^{\vep,n}(t,.)-\vpo\|_1>\a\}$,
where $u^{\vep,n}$ is here the solution of equation (\ref{1}), but with an additive noise
that becomes stationary in space as $n$ goes to infinity (see \cite{dBG} for a precise statement)
then there
exists a constant $C(\alpha,c_0)$ which depends on $\alpha$ and
$c_0$ but not on $T$ such that
\begin{equation}
\label{ssmod}
\underline{\lim}_{n\rightarrow\infty}\underline{\lim}_{\vep\to0}\vep^2\log\mathbb{P}
\left(\tilde{\tau}_{\alpha}^{n,\vep} \le
T\right)\ge-\frac{C(\alpha,c_0)}{T^3}.
\end{equation}
It is not clear that (\ref{ssmod}) is still true in the present multiplicative case,
because the proof involves a controlability problem with a potential which -- up to now --
is open.

Note also that the decomposition given in Theorem \ref{t2}
is not unique, and is determined by the choice of specific orthogonality conditions
(see the proof below). In particular, contrary to the additive case, we will
be able here to investigate the asymptotic behavior in time of the limit
process by choosing one particular decomposition of the form given in Theorem \ref{t2}.
This is the object of Section 3.3.

\noindent
{\em Proof of Theorem \ref{t2}}
The proof follows closely the proof of Theorem 2.1 in \cite{dBD1} and we refer to \cite{dBD1}
for more details. The parameters $\xet$ and $\cet$ are obtained thanks to the use of the implicit
function Theorem. These are then local semi-martingales defined as long as
$|\cet-c_0|<\a$ and $\|\ue(t,.+\xet)-\vpo\|_1<\a$, and setting
$$ \vep \ee(t)=\ue(t,.+\xet)-\pet,$$
one has for each $\vep>0$, almost surely,
\begin{equation}
\label{ortho}
(\ee,\vpo)=(\ee, \p_x \vpo)=0.
\end{equation}
In order to estimate the exit time 
$$
\tea=\inf\{t\ge 0, \; |\cet-c_0|>\a \; \mathrm{or} \; \|\vep \ee (t) \|_1 >\a\},
$$
we make use , as in \cite{dBD1}, of the functional defined for $u \in H^1(\R)$,
\begin{equation}
\label{lyap}
Q_{c_0}(u):= H(u)+c_0 m(u)
\end{equation}
where $H$ and $m$ are defined respectively in (\ref{masse}) and (\ref{hamiltonien}).
Note that $\vpo$ is a critical point of $Q_{c_0}$. We denote by $L_{c_0}$ the linearized operator
around $\vpo$, that is
\begin{equation}
\label{linearise}
L_{c_0}=-\p_x^2+c_0-2\vpo.
\end{equation}
The next lemma, which is proved with the use of the It\^o Formula, using the same regularization procedure
as in \cite{dBD0}, gives the evolution of $H$ and $m$
for the solution $\ue$ of (\ref{1}) with $\ue(0)=\vpo$~:

\begin{Lemma}
\label{l2.1}
For any stopping time $\tau <+\infty$ a.s, one has
$$
m(\ue(\t))= m(\vpo) -\vep \int_0^{\t} ((\ue)(s), dW(s)) + \vep^2 |k|^2_{L^2} \int_0^{\t} m(\ue(s))ds
$$
and
\begin{eqnarray}
H(\ue(\t))&=& H(\vpo)+\vep \int_0^{\t} (\p_x\ue, \p_x(\ue dW(s)))
-\frac{\vep}{2} \int_0^{\t} ((\ue)^3,dW(s))\\
& &+\frac{\vep^2}{2} \int_0^{\t} \left\{ |k|_{L^2}^2 |\p_x \ue|_{L^2}^2 
+ |k'|_{L^2}^2 |\ue|_{L^2}^2\right\}ds \\
& & - \frac{\vep^2}{2} \Sum_k \int_0^{\tau} \int_{\R} (\ue)^3 |\phi e_k|^2 dx ds.
\end{eqnarray}
\end{Lemma}

Consider $\nu>0$ such that  $(Q''_{c_0}(\vpo)v,v)\ge \nu \|v\|_1^2$ for any $v \in H^1$
satisfying $(v,\vpo)=(v,\p_x \vpo)=0$. The existence of such a constant is a classical result
(see \cite{Ben} or \cite{BSS}). Then it is easy to show (see \cite{dBD1}) that there is a constant
$C(\a_0)>0$ such that for any $t<\tea$,
\begin{equation}
\label{minq}
Q_{c_0}(\ue(t,.+\xet))-Q_{c_0}(\pet)\ge \frac{\nu}{4} \|\vep \ee(t)\|_1^2 -C|\cet-c_0|^2.
\end{equation}
Now, if $\t=\tea\wedge t$, then by (\ref{minq}), the
translation invariance of $Q_{c_0}$, and Lemma \ref{l2.1}
\begin{equation}
\label{majeta}
\begin{array}{rcl}
\|\vep \ee(\t)\|_1^2 &\le & \Frac{4}{\nu} \left[ Q_{c_0}(\vpo) -Q_{c_0}(\peto)\right]
+\vep \Int_0^{\t } (\p_x\ue(s),\p_x(\ue dW(s))) \\
& & -\Frac{\vep}{2} \Int_0^{\t} ((\ue)^3(s),dW(s)) +\Frac{\vep^2}{2} \Int_0^{\t}
(|k|_{L^2}^2 |\p_x \ue|_{L^2}^2 +|k'|_{L^2}^2 |\ue|_{L^2}^2) ds \\
& & -\Frac{\vep^2}{2} \Sum_k \int_0^{\t} \int_{\R} (\ue)^3(s) |\phi e_k|^2 dx ds -c_0 \vep \Int_0^{\t} ((\ue)^2,dW(s)) \\
& & +c_0\vep^2 |k|_{L^2}^2 \Int_0^{\t} m(\ue(s))ds +C |\ceto-c_0|^2.
\end{array}
\end{equation}
The term $|\ceto-c_0|$ is then estimated thanks to the orthogonality condition $(\ee,\vpo)=0$
and the evolution of $m(\ue(\t))$ given in Lemma \ref{l2.1}; one obtains, for some constants
$\mu>0$ and $C>0$, depending only on $c_0$ and $\a_0$ (with $\a\le \a_0$)
\begin{eqnarray*}
\mu |\ceto-c_0| &\le & \left| |\vpo|_{L^2}^2 - |\peto|_{L^2}^2\right| \\
& \le & |\vep \ee(\t)|_{L^2}^2 +C\a|\ceto-c_0|+2\vep \left| \int_0^{\t} ((\ue)^2,dW(s)) \right|\\
& & +2 \vep^2 |k|_{L^2}^2 \int_0^{\t} |\ue(s)|_{L^2}^2 ds.
\end{eqnarray*}
Hence, choosing $\a_0$ sufficiently small one gets
\begin{equation}
\label{majc}
\begin{array}{rcl}
|\ceto-c_0|^2 & \le & C \Big[ |\vep \ee(\t)|_{L^2}^4 +4\vep^2 \Big| \Int_0^{\t}
((\ue)^2,dW(s)) \Big|^2\\
& & +4 \vep^4 |k|_{L^2}^4 \Big( \Int_0^{\t} |\ue(s)|_{L^2}^2ds\Big)^2 \Big]
\end{array}
\end{equation}
which, once inserted into (\ref{majeta}) leads to
\begin{eqnarray*}
\|\vep \ee(\t)\|_1^2 & \le & C \Big[ |\vep \ee(\t)|_{L^2}^4 +\vep \Big| \int_0^{\t}
(\p_x \ue,\p_x (\ue dW(s)))\Big| \\
& & +\vep \Big| \int_0^{\t} ((\ue)^3,dW(s))\Big| +c_0 \vep \Big| \int_0^{\t} ((\ue)^2,dW(s)) \Big|\\
& & +4 \vep^2 \Big| \Int_0^{\t} ((\ue)^2, dW(s))\Big|^2 +\vep^2 \|k\|_1^2 \int_0^{\t} \|\ue(s)\|_1^2 ds \\
& &+ \vep^2 |k|^2_{L^2} \int_0^{\t} \|\ue(s)\|_1^3 ds +\vep^4 |k|_{L^2}^4 \Big( \int_0^{\t}
|\ue(s)|_{L^2}^2ds\Big)^2 \Big].
\end{eqnarray*} 
With this estimate in hand, together with (\ref{majc}), the  conclusion of Theorem \ref{t2}
follows with the same arguments as in the proof of Proposition 3.1 in \cite{dBG}.
These arguments rely on classical exponential tail estimates for stochastic integrals,
after noticing that $\|\ue(s)\|_1\le C$, a.s. for $s\in [0,\tea\wedge T]$ and $\a\le\a_0$,
so that the quadratic variation of each of the integrals involved in the above estimates are
bounded above by $CT$.
\hfill
$\square$

\section{A central limit theorem}
This section is devoted to the proof of the next theorem:

\begin{Theorem}
\label{t3}
Under the assumptions of Theorem \ref{t2}, let $\a<\a_0$ be fixed.
Then we can find ${\tilde c}^{\vep}(t)$ and ${\tilde x}^{\vep}(t)$ 
satisfying the conclusion of Theorem \ref{t2}
such that if ${\tilde \e}^{\vep}$ is defined as in Theorem \ref{t2}, for any $T>0$, the
process $({\tilde \e}^{\vep}(t))_{t \in [0,T]}$ converges in $L^2(\O; L^{\infty}(0,\tea\wedge
T;L^2(\R)))$ to a Gaussian process $\tilde \e$ satisfying
the additive linear equation
\begin{equation}
\label{eql}
d\tilde\e= \p_x L_{c_0}\tilde \e \, dt + \tilde Q\vpo d\tilde W,
\end{equation}
with $\tilde \e(0)=0$, where $ \tilde W$ is the Wiener process with covariance
$\phi\phi^*$ given by $\tilde W=\T_{c_0 t} W$, and $\tilde Q$ is a
projection operator.
Moreover, for $a>0$ sufficiently small compared to $c_0$, the process
$w(t,x)=e^{ax}\tilde \e(t,x)$ is a well defined $H^1$ valued process, of
Ornstein-Uhlenbeck type, which converges in law to an $H^1$-valued Gaussian
random variable as $t$ goes to infinity.
\end{Theorem}

The conclusion of Theorem \ref{t3} will be obtained in three steps.
The first step consists in estimating the modulation parameters obtained in
Theorem \ref{t2}, in terms of $\ee$, using the equations for those parameters;
then the convergence of $\ee$ as $\vep$ tends to zero is proved, and finally
in the third step, a slight change in the modulation parameters is performed,
in order that the limit process $\e$ may be written as an Ornstein-Uhlenbeck
process.

From now on, we assume that $\a$ is fixed and sufficiently small, so that the conclusion
of Theorem \ref{t2} holds, and we denote $\tea$ by $\te$.

\subsection{Modulation equations}

Since we know that the modulation parameters $\xet$ and $\cet$ are semi-martingale processes
adapted to the filtration generated by $(W(t))_{t\ge 0}$, we may a priori write the stochastic
evolution equations for those parameters in the form
\begin{equation}
\label{modeq}
\left\{ 
\begin{array}{l}
d\xe = \ce dt +\vep \ye dt +\vep (\ze, dW)\\
d\ce = \vep \ae dt +\vep (\be, dW)
\end{array}
\right.
\end{equation}
where $\ye$ and $\ae$ are real valued adapted processes with a.s. locally integrable paths
on $[0,\te)$,
and $\be$, $\ze$ are predictable processes with paths a.s. in $L^2_{loc}(0,\te; L^2(\R))$.
We then proceed as in \cite{dBD1} : the It\^o-Wentzell Formula applied to $\ue(t,x+\xet)$,
together with equation (\ref{1}) for $\ue$ and the first equation of (\ref{modeq}) for $\xe$ give a stochastic
evolution equation for $\ue(t, x+\xe)$.
On the other hand, the standard It\^o Formula together with the second equation of (\ref{modeq}) for $\ce$ 
give an equation for the evolution of $\pet$.
Replacing then $\pet+\vep \ee(t,x)$ for $\ue(t,x+\xet)$ in the first equation leads to the 
following stochastic equation for the evolution of $\ee(t)$ :
\begin{equation}
\label{eqetae}
\renewcommand{\arraystretch}{1.7}
\begin{array}{rcl}
d\ee & = & \p_x L_{c_0} \ee dt +(\ye\p_x \pe -\ae \p_c \pe) dt- \p_x ( (\pe-\vpo)\ee) dt \\
& & + (\ce-c_0+\vep \ye)\p_x \ee dt -\frac{\vep}{2} \p_x ((\ee)^2)dt +\pe \T_{\xe}dW \\
& & + \p_x\pe (\ze,dW) -\p_c \pe (\be,dW) +\vep \ee \T_{\xe}dW +\vep \p_x \ee (\ze, dW) \\
& & +\frac{\vep}{2} \p_x^2 \pe |\phi^*\ze|_{L^2}^2 dt -\frac{\vep}{2} \p_c^2 \pe |\phi^* \be|_{L^2}^2 dt+\vep \Sum_{l\in \N} \p_x(\pe \T_{\xe}\phi e_l) (\ze,\phi e_l) dt \\
& & +\frac12 \vep^2 \p_x^2 \ee |\phi^* \ze|_{L^2}^2 dt +\vep^2 \Sum_{l \in \N} \p_x (\ee \T_{\xe}\phi e_l)(\ze, \phi e_l) dt
\end{array}
\end{equation}
where $L_{c_0}$ is defined in (\ref{linearise}).
Now, taking the $L^2$- inner product of equation (\ref{eqetae}) with $\vpo$, on the one hand, 
and with $\p_x \vpo$ on the other hand, then using the orthogonality conditions (\ref{ortho})
and the fact that $L_{c_0} \p_x \vpo =0$, and finally identifying the drift parts and the martingale
parts of each of the resulting equations lead to the same kind of system that we previously obtained
in \cite{dBD1}; namely, setting
$$
Y^{\vep}(t)= \begin{pmatrix} \ye(t)\\ \ae(t) \end{pmatrix}
\quad \mathrm{and} \quad Z^{\vep}_l(t) =\begin{pmatrix} (\ze, \phi e_l)\\(\be, \phi e_l)\end{pmatrix}
$$
then one gets for the drift parts
\begin{equation}
\label{eqdrift}
A^{\vep}(t) Y^{\vep}(t) = G^{\vep}(t)
\end{equation}
where
\begin{equation}
\label{eqA}
A^{\vep}(t)= \begin{pmatrix} (\p_x\pe+\vep \p_x \ee, \p_x \vpo) & -(\p_c \pe,\p_x\vpo) \\
                            -(\p_x \pe, \vpo) & (\p_c \pe, \vpo)
            \end{pmatrix}
\end{equation}
and 
$$
G^{\vep}(t)=\begin{pmatrix}
            G^{\vep}_1(t) \\ G^{\vep}_2(t)
           \end{pmatrix},
$$
with
\begin{equation}
\label{eqG1}
\renewcommand{\arraystretch}{1.7}
\begin{array}{rcl}
G^{\vep}_1(t) &=& (\ee, L_{c_0}\p_x^2 \vpo)+(\ce-c_0)(\ee, \p_x^2 \vpo) +\frac{\vep}{2} (\p_x(\ee)^2, \p_x \vpo)\\
& &  +(\p_x((\pe-\vpo)\ee), \p_x\vpo) -\frac{\vep}{2} (\p_x^2 \pe, \p_x \vpo)|\phi^* \ze|_{L^2}^2 \\
& & +\frac{\vep}{2} (\p_c^2 \pe, \p_x \vpo) |\phi^* \be|_{L^2}^2 -\vep \Sum_{l\in \N} (\ze, \phi e_l)(\p_x(\pe \T_{\xe}\phi e_l), \p_x \vpo) \\
& & +\frac12 \vep^2 (\ee, \p_x^3 \vpo) |\phi^* \ze|_{L^2}^2 -\vep^2 \Sum_{l \in \N} (\p_x(\ee \T_{\xe}\phi e_l), \p_x \vpo) (\ze, \phi e_l)
\end{array}
\end{equation}
and
\begin{equation}
\label{eqG2}
\renewcommand{\arraystretch}{1.7}
\begin{array}{rcl}
G^{\vep}_2(t)&=& -\frac{\vep}{2} (\p_x(\ee)^2, \vpo)-(\p_x((\pe-\vpo)\ee), \vpo) +\frac{\vep}{2} (\p_x^2 \pe,\vpo)|\phi^* \ze|_{L^2}^2\\ 
& & -\frac{\vep}{2} (\p_c^2 \pe, \vpo)|\phi^* \be|_{L^2}^2 +\vep \Sum (\ze,\phi e_l) (\p_x(\pe\T_{\xe}\phi e_l) ,\vpo) \\
& & +\frac{\vep^2}{2} (\ee,\p_x^2 \vpo)|\phi^* \ze|_{L^2}^2 +\vep^2 \Sum_{l\in \N} (\p_x (\ee \T_{\xe} \phi e_l), \vpo) (\ze, \phi e_l);
\end{array}
\end{equation}
note that $A^{\vep}(t)=A_0+O(|\ce-c_0|+\|\vep \ee\|_1)$, a.s. for $t\le \te$
with
$$
A_0= \begin{pmatrix}
     |\p_x \vpo|^2_{L^2} & 0 \\ 0  & (\vpo, \p_c \vpo)
    \end{pmatrix}
$$
and $O(|\ce-c_0|+\|\ee\|_1)$ is uniform in $\vep, t$ and $\o$ as long as $t\le \te$.
Concerning the martingale parts, one gets the equation
\begin{equation}
\label{eqmart}
A^{\vep}(t) Z^{\vep}_l(t) = F^{\vep}_l(t), \quad \forall l \in \N
\end{equation}
with 
\begin{equation}
\label{eqF}
F ^{\vep}(t) =\begin{pmatrix}
-((\pe+\vep \ee)\T_{\xe}\phi e_l, \p_x \vpo) \\ ((\pe+\vep\ee) \T_{\xe} \phi e_l, \vpo).
             \end{pmatrix}
\end{equation}
\begin{Proposition}
\label{p1}
Under the above assumptions, there is a constant $\a_1>0$, such that if $\a \le \a_1$, then
\begin{equation}
\label{majbz}
|\phi^* \ze(t)|_{L^2} +|\phi^* \be |_{L^2} \le C_1|k|_{L^2}, \quad \mathrm{a.s.}\; \mathrm{for} \; t\le \te
\end{equation}
and
\begin{equation}
\label{majay}
|\ae(t)| +|\ye(t)| \le C_2 |\ee(t)|_{L^2} +\vep C_3, \quad \mathrm{a.s.}\; \mathrm{for} \; t\le \te
\end{equation}
for some constants $C_1$, $C_2$, $C_3$, depending only  on $\a$ and $c_0$, and for any $\vep \le \vep_0$.
\end{Proposition}

\noindent
{\em Proof}
The proof is exactly the same as the proof of  Corollary 4.3 in \cite{dBD1}, once noticed that, a.s. for $t\le \te$,
\begin{eqnarray*}
\Sum_{l\in \N} |F^{\vep}_l(t)|^2 & \le & C \Sum_{l \in \N} |(\pe+\vep \ee) \T_{\xe}\phi e_l|_{L^2}^2 \\
& \le& C \Sum_{l} \int_{\R} (\pe+\vep \ee)^2 (x) [ (\T_{\xe}k)* e_l]^2(x) dx\\
& \le & \int_{\R} (\pe+\vep \ee)^2 (x)\Sum_l (\T_{\xe} k(x-.), e_l)^2dx\\
& \le & C \int_{\R} (\pe+\vep \ee)^2 (x) |\T_{\xe} k(x-.)|_{L^2}^2 dx\\
& \le & C|k|_{L^2}^2 |\pe+\vep \ee|_{L^2}^2 \le C|k|_{L^2}^2
\end{eqnarray*}
where we have used the Parseval equality in the fourth line.
\hfill 
$\square$

\subsection{Convergence of $\ee$}

Let us first assume that $\ee$ has a limit as $\vep$ goes to zero, and take formally the 
limit as $\vep$ goes to zero in the preceding equations. Then, as was noticed above,
$$
\lim_{\vep \to 0} A^{\vep} =A_0= \begin{pmatrix}
                                |\p_x \vpo|_{L^2}^2 & 0 \\ 0 & (\vpo, \p_c \vpo) 
                               \end{pmatrix}
$$
hence 
\begin{equation}
\label{limz}
\lim_{\vep \to 0} \phi^* \ze = -\frac{1}{|\p_x\vpo|_{L^2}^2} (\T_{c_0 t}\phi)^* (\vpo \p_x \vpo) :=z
\end{equation}
\begin{equation}
\label{limb}
\lim_{\vep \to 0} \phi^* \be= \frac{1}{(\vpo,\p_c \vpo)} (\T_{c_0 t}\phi^*) (\vp_{c_0}^2) :=b
\end{equation}
\begin{equation}
\label{limy}
\lim_{\vep \to 0} y^{\vep}= \frac{1}{|\p_x \vpo|_{L^2}^2} (\e, L_{c_0} \p_x^2 \vpo):=y
\end{equation}
and
\begin{equation}
\label{lima}
\lim_{\vep \to 0} \ae=0.
\end{equation}
Moreover, formally, $\e$ satisfies the equation
\begin{equation}
\label{eqeta}
\begin{array}{rcl}
d\e & = & \p_x L_{c_0} \e dt + \frac{1}{|\p_x \vpo|_{L^2}^2} (\e, L_{c_0}\p_x^2 \vpo) \p_x \vpo dt \\
& & + \vpo \T_{c_0 t} dW -\frac{1}{2 |\p_x \vpo|_{L^2}^2} (\p_x(\vp_{c_0}^2), \T_{c_0 t} dW)\p_x \vpo\\
& & - \frac{1}{(\vpo, \p_c \vpo) } (\vp_{c_0}^2, \T_{c_0 t} dW) \p_c \vpo.
\end{array}
\end{equation}
It is easy to show that (\ref{eqeta}) has a unique adapted solution $\e$ with paths a.s. in $C(\R^+, H^1)$ satisfying
$\e(0)=0$. Moreover using the fact that $(\p_c \vpo, \p_x \vpo)=0$, one easily gets from the above equation that 
$ (\e, \vpo)=(\e,\p_x \vpo)=0$, $\forall t>0$.

Next, we make use of the following lemmas, whose proofs are obtained in the same way as the corresponding Lemmas in
\cite{dBD1}.

\begin{Lemma}
\label{l2}
Let $\e$ be the solution of (\ref{eqeta}) with $\e(0)=0$.
Then, for any $T>0$, there is a constant $C$ depending only on $c_0$, $T$ and $\|k\|_1$ such that 
$$
\E\left( \|\e(t)\|_1^4\right) \le C, \; \forall t\le T.
$$
\end{Lemma}




\begin{Lemma}
\label{l3}
Let $\ee$ be the solution of (\ref{eqetae}), defined for $t \in [0,\te[$, obtained
thanks to the modulation procedure of Section 2.
Then, for any $T>0$, 
$$
\E\Big( \sup_{t\le\te\wedge T} |\ee(t) |_{L^2}^4 \Big) \le C(T,\a, c_0,\|k\|_1).
$$
\end{Lemma}


The above lemmas show that 
\begin{equation}
\label{cvgencec}
\forall T>0, \; \forall q\ge 2, \; \lim_{\vep \to 0} \E \Big( \sup_{t\le T\wedge \te} |\cet -c_0|^q\Big)=0.
\end{equation}
Indeed, the expression of $\cet-c_0$ given by (\ref{modeq}) together with (\ref{majbz}) and (\ref{majay}) imply
easily 
$$
\E \Big( \sup_{t\le T\wedge \te} |\cet -c_0|^2\Big) \le C\vep^2 [1+\E \int_0^{T\wedge \te} |\ee(s)|_{L^2}^2ds]
$$
with $C=C(\a, c_0,T,\|k\|_1)$. Then, (\ref{cvgencec}) is deduced form Lemma \ref{l3} for $q=2$, and follows 
for all other values of $q$ from the uniform boundedness of $|\cet-c_0|$ on $[0,T\wedge\te]$.
Note that an immediate consequence of (\ref{cvgencec}) is the fact that
\begin{equation}
\forall T>0, \; \forall q\ge 2, \;  \lim_{\vep \to 0} \E \Big( \sup_{t\le T\wedge \te} \|\pet-\vpo\|_2^2\Big)=0.
\end{equation}

We will finally need the next lemma.

\begin{Lemma}
\label{l6}
For any $T>0$, and any $q\ge 1$,
$$
\lim_{\vep \to 0} \E \Big( \sup_{t\le T\wedge \te} \Big( \sum_{l\in \N}|Z^{\vep}_l(t) -Z_l(t)|^2\Big)^q \Big)=0
$$
where we have set for $l \in \N$
$$
Z_l(t)=\begin{pmatrix} 
       (z,\phi e_l) \\ (b,\phi e_l)
      \end{pmatrix},
$$
$z$ and $b$ being given by (\ref{limz}) and (\ref{limb}), respectively.
\end{Lemma}

\noindent
{\em Proof}
Here again, it is sufficient to consider the case $q=1$. We recall that 
$Z^{\vep}_l$ satisfies equation (\ref{eqmart}).
First, it is clear that 
$$
\lim_{\vep \to 0} \E \Big( \sup_{t\le T\wedge \te} \|(A^{\vep}(t))^{-1} -(A_0(t))^{-1}\|^{2q}\Big) =0, \; \forall q\ge 1.
$$
On the other hand, in view of (\ref{eqF}), denoting $F^0_l(t)$ the formal limit
of $F^{\vep}_l(t)$, one has
$$
\begin{array}{l}
\E \Big( \Sup_{t\le T\wedge\te} \sum_l |F^{\vep}_l(t)-F^0_l(t)|^2\Big) \\
\le C \E \Big( \Sup_{t\le T\wedge\te} \sum_l |\p_x \vpo(\T_{\xe}\phi-\T_{c_0 t} \phi)e_l|_{L^2}^2\Big)
+C \E \Big( \Sup_{t\le T\wedge\te} \|\pet-\vpo\|_1^2\Big)
\end{array}
$$
and
$$
\begin{array}{l}
\E \Big( \Sup_{t\le T\wedge\te} \Sum_l |\p_x \vpo(\T_{\xe}\phi-\T_{c_0 t} \phi)e_l|_{L^2}^2\Big)\\
\le \|\vpo\|_1^2 \E \Big( \Sup_{t\le T\wedge\te} |k(.+\xet -c_0 t)-k|_{L^2}^2\Big).
\end{array}
$$
Then, the It\^o Formula applied to the function 
$$
\K^{\vep}(t,x)=(k(x+\xet-c_0 t)-k(x))^2
$$
using equation (\ref{modeq}) for $d\xet$, together with (\ref{majbz}), (\ref{majay}),
and (\ref{cvgencec}) lead to the conclusion
of Lemma \ref{l6}.
\hfill $\square$

Now, in order to prove that
\begin{equation}
\label{cvgenceeta} 
\lim_{\vep \to 0} \E \Big( \sup_{t\le T\wedge\te} |\ee(t)-\e(t)|_{L^2}^2\Big) =0,
\end{equation}
where $\e$ is the solution of (\ref{eqeta}) with $\e(0)=0$, it suffices to set $v^{\vep}=\ee-\e$,
to deduce from (\ref{eqeta}) and (\ref{eqetae}) the equation for $dv^{\vep}$ and to apply the
It\^o Formula to get the evolution of $|v^{\vep}|_{L^2}^2$.
We do not give the details of those tedious, but easy computations.
Finally, the use of the following estimates :
$$
\begin{array}{l}
\vep |(\ve, \p_x ((\ee)^2))|=\vep |(\p_x \e,(\ee)^2)|\le \vep \|\e\|_1 |\ee|_{L^4}^2 \\
\le C \vep \|\e\|_1 |\ee|_{L^2}^{3/2}|\p_x\ee|_{L^2}^{1/2}
\le C \sqrt{\vep}\|\e\|_1 |\ee|_{L^2}^{3/2}
\end{array}
$$
on the one hand, and
$$
|y^{\vep} -y| +|a^{\vep}|\le  C (|\ve|_{L^2} + |\ce -c_0||\ee|_{L^2} +\vep |\ee|_{L^2}^2 
+|\ee|_{L^2} \|\pe -\vpo\|_1 +\vep)
$$
which is obtained as in the proof of Lemma \ref{l6} on the other hand, together with Lemma \ref{l2}
to \ref{l6} allow to get the conclusion, that is the convergence of $\ee$ to $\e$ in
$L^2(\O,L^{\infty}(0,\te \wedge T;L^2(\R)))$.
\hfill
$\square$

\subsection{Complements on the limit equation}

First of all, we note that the modulation equations may be written at order one in $\vep$
as
$$
\left\{ 
\begin{array}{l}
d\xe = c_0 dt +\vep y dt +\vep W_1 dt +\vep dW_2 +o(\vep)\\
d\ce = \vep dW_1 +o(\vep)
\end{array}
\right.
$$
where
$$
y=|\p_x \vpo|_{L^2}^{-2} (\e,L_{c_0}\p_x^2 \vpo),
$$
$$
W_1(t)=(\vpo,\p_c \vpo)^{-1} (\vp_{c_0}^2,\tilde W(t))
$$
and
$$
W_2(t)=-\frac12 |\p_x\vpo|_{L^2}^{-2} (\p_x (\vpo^2),\tilde W(t)).
$$

Note that $W_1$ and $W_2$ are real valued Brownian motions, which are independent
since
$$
\E(W_1(t)W_2(s)) =-\frac12 |\p_x\vpo|_{L^2}^{-2} (\vpo,\p_c\vpo)^{-1} (\phi^*(\p_x(\vp_{c_0}^2)),\phi^*(\vp_{c_0}^2))(t\wedge s)
= 0
$$
because the operator $\phi^*$ commutes with spatial derivation.

Now, we want to investigate the asymptotic behavior in time of the process $\eta$.
However, in the present form, the process $\e$ does not converge in law as $t$ goes to infinity;
this is due to the fact that the preceding modulation does not exactly correspond to the projection
of the solution $\ue$ on the (two-dimensional) center manifold, in which case the remaining term would
belong to the stable manifold around the soliton trajectory.
We now show that by slightly changing the modulation parameters, we can get a new decomposition of
the solution $\ue$ which is defined on the same time interval as before, 
but which fits with the preceding requirements.
For that purpose, we first need to recall a few facts from \cite{PW}.

The generalized nullspace of the operator $\p_x L_{c_0}$ (that is the operator arising in
the linearized evolution equation in the soliton reference frame)
is spanned by the functions $\p_x \vpo$ and $\p_c \vpo$, with the equality
$$
\p_x L_{c_0} \p_c \vpo= -\p_x \vpo
$$
and there are constants $\theta_1$ and $\theta_2$ (with $\theta_1=(\vpo,\p_c\vpo)$) such that
if we set
$$
\tilde g_1(x)=-\theta_1\int_{-\infty}^x \p_c\vpo(y)dy +\theta_2 \vpo \quad \mathrm{and} \quad \tilde g_2(x)=\theta_1\vpo
$$
then the generalized nullspace of $-L_{c_0}\p_x$ is spanned by $\tilde g_1$ and $\tilde g_2$ and 
$$
(\tilde g_1,\p_x \vpo)=1, \; (\tilde g_1, \p_c \vpo)=0, \; (\tilde g_2, \p_x \vpo)=0, \; (\tilde g_2, \p_c \vpo)=1.
$$
We also set, for $a>0$,
$$
f_1^a(x)=e^{ax}\p_x \vpo, \; f_2^a(x)=e^{ax} \p_c \vpo, \; g_1^a(x)=e^{-ax}\tilde g_1(x),  \; g_2^a(x)=e^{-ax}\tilde g_2(x),
$$
so that $(f_i^a,g_j^a)=\delta_{ij}$. Then the operator $A_a$ defined for 
$a>0$ by $A_a=e^{ax}\p_x L_{c_0}e^{-ax}$ has a well defined generalized nullspace
spanned by $f_1^a, f_2^a$ and the spectral projection on this nullspace is
given by $Pw=\sum_{k=1}^2 (w,g_k^a)f_k^a$ where $w=e^{ax} v$, and $v$ is an $L^2$
function. Moreover, if $Q=I-P$, then $Q$ is the spectral projection on the stable
manifold of $A_a$, and under the condition $0<a<\sqrt{c_0/3}$, there are constants 
$C>0$ and $b>0$ such that
\begin{equation}
\label{decexp}
\|e^{A_a t}Qw\|_1 \leq Ce^{-bt}\|w\|_1, \; \forall t>0, \; \forall w \in H^1,
\end{equation}
where $e^{A_a t}$ is the $C^0$-semi-group generated by $A_a$ (see Theorem 4.2 in \cite{PW}).

Now, let $\e$ be the solution of (\ref{eqeta}) with $\e(0)=0$, and consider $w(t,x)=e^{ax}\e(t,x)$.
Note that the orthogonality condition $(\e,\vpo)=0$ implies $(w,g_2^a)=0$, so that
$Pw=\lambda (t) f_1^a$ with $\lambda(t)=(w(t),g_1^a)$ a real valued stochastic process
whose evolution is given by
\begin{equation}
\label{eqlam}
\begin{array}{rcl}
\lambda(t)&=& \Int_0^t |\p_x \vpo|_{L^2}^{-2} (\e(s),L_{c_0}\p_x^2\vpo)ds
-\Int_0^t |\p_x \vpo|_{L^2}^{-2}(\vpo\p_x\vpo, d\tilde W(s)) \\ & &  +\Int_0^t (e^{ax}\vpo d\tilde W(s),g_1^a)
\end{array}
\end{equation}
where we have used (\ref{eqeta}) and the fact that $A_a Pw=0$ and $\lambda(0)=0$. Hence, $\lambda(t)$ is bounded
in $L^4(\O;L^{\infty}(0,T\wedge \te))$ by Lemma \ref{l2}.
Let us set ${\tilde x}^{\vep}(t)=\xet-\vep \lambda(t)$ for $t\in [0,\te[$.
Then
\begin{equation}
\label{newmod}
\ue(t,x+{\tilde x}^{\vep}(t))=\pet(x)+\vep {\tilde \e}^{\vep}(t,x)
\end{equation}
with
$$
{\tilde \e}^{\vep}(t,x)=\frac{1}{\vep} (\pet(x-\vep\lambda(t))-\pet(x))+\ee(t,x-\vep\lambda(t)).
$$
Note that, a.s. for $t\le \te$ :
$$
|\pet(.-\vep \lambda(t))-\pet-\vep \lambda(t) \p_x \pet|_{L^2} 
\le \vep^2 \lambda^2(t) C(c_0,\a).
$$
Hence, it follows from Lemma \ref{l2}, \ref{l3} and the above bound on $\lambda$
that
\begin{equation}
\label{limetatilde}
\lim_{\vep \to 0} \E \big( \sup_{t\le T\wedge \te} |{\tilde \e}^{\vep}(t)-\tilde \e(t)|_{L^2}^2 \big)=0
\end{equation}
with $\tilde \e(t)=\e(t)-\lambda(t)\p_x \vpo$.
So now, with this new decomposition, we clearly have, setting $\tilde w(t,x)=e^{ax}\tilde \e(t,x)$ :
$$
P\tilde w=0, \quad Q\tilde w=Qw.
$$
Also, if $w_2=Qw$, then the equation (\ref{eqeta}) implies
\begin{equation}
\label{eqw2}
dw_2=A_aw_2 dt +Qe^{ax}\vpo d\tilde W
\end{equation}
hence
$$
w_2(t)=\int_0^t e^{A_a(t-\s)}Q[e^{ax}\vpo d\tilde W(\s)];
$$
the trace of the covariance operator of the Gaussian process $w_2$ in $H^1$ may be easily computed and
estimated thanks to (\ref{decexp}) as
$$
\int_0^t \sum_l \|e^{A_a \s} Q e^{ax}\vpo\phi e_l\|_1^2 d\s
\le C \big(\int_0^t e^{-b\s}d\s\big) \sum_l \|e^{ax}\vpo\phi e_l\|_1^2 d\s
\le C \|k\|_1^2 \|e^{ax}\vpo\|_1^2.
$$
Moreover, this covariance operator converges as $t$ goes to infinity and it follows that $w_2$ converges
in law in $H^1$ to a Gaussian random variable.
The end of the statement of Theorem \ref{t3} follows, setting $\tilde Qv=e^{-ax}Qe^{ax}v$.
\hfill
$\square$

\section{A remark on the soliton diffusion}

Let us go back to the stochastic evolution equations for the new modulation parameters, that we may write as
\begin{equation}
\label{newmodeq}
\left\{ 
\begin{array}{l}
d{\tilde x}^{\vep} = c_0 dt +\vep B_1 dt +\vep dB_2 +o(\vep) \\
d\ce = \vep dB_1 +o(\vep)
\end{array}
\right.
\end{equation}
with $B_1=W_1$ and $B_2=-(e^{ax}\vpo \tilde W(t),g_1^a)=-(\tilde W(t), \vpo \tilde g_1)$.
Note that $B_1$ and $B_2$ are now correlated Brownian motions. We denote by 
$$\s=(\s_{ij})_{i,j}=\mathrm{cov}(B_1,B_2).
$$
If we keep only the order one terms in $\vep$ i.e. we consider the solution $(X^{\vep}(t), C^{\vep}(t))$ of the system of SDEs
$$
\left\{
\begin{array}{l}
dX^{\vep}= c_0dt +\vep B_1dt +\vep dB_2\\
dC^{\vep}=\vep dB_1,
\end{array}
\right.
$$
then $(X^{\vep}(t)-c_0t, C^{\vep}(t)-c_0)$ is a centered Gaussian vector, and it is easy to compute its covariance
matrix. Let us denote by $\mu^{\vep}_t$ the law of $(X^{\vep}(t)-c_0t, C^{\vep}(t)-c_0)$; we may compute
\begin{equation}
\label{maxesp}
\begin{array}{l}
\displaystyle \max_{x\in \R} \E \Big( \vp_{C^{\vep}(t)}(x-X^{\vep}(t)) \Big) \\
= \displaystyle \max_{x\in \R } \Int \!\! \Int \vp_{c+c_0}(x-c_0 t-y)
\mu^{\vep}_t (dy,dc)\\
= \displaystyle \max_{x\in \R } \Frac{1}{(\det \Sigma)^{1/2}}\Int\!\! \Int  \vp_{c+c_0} (x-c_0 t -y) \exp \Big(-\frac12 \Sigma^{-1}\begin{pmatrix}
                                                                                     c\\y
                                                                                    \end{pmatrix}.
\begin{pmatrix}
c\\y
\end{pmatrix}\Big)dc dy
\end{array}
\end{equation}
where $\Sigma$ is the covariance matrix of $(X^{\vep}(t)-c_0t, C^{\vep}(t)-c_0)$, given by
$$
\Sigma = \vep^2 \begin{pmatrix}
                \s_{11} t & \s_{12} t +\s_{11} \frac{t^2}{2} \\
		\s_{12} t +\s_{11} \frac{t^2}{2} & \s_{22} t +\s_{12} t^2 +\s_{11} \frac{t^3}{3}
               \end{pmatrix}.
$$
It is not difficult to see that
$$
\exp \Big( -\frac12 \Sigma^{-1} \begin{pmatrix}
                                c\\y
                               \end{pmatrix}.
\begin{pmatrix}
c\\y
\end{pmatrix}\Big)
\le \exp\Big( -\frac12 \frac{\vep^2}{\det \Sigma} \big( \s_{11} \frac{t^3}{12} +(\s_{22}-\frac{\s_{12}^2}{\s_{11}}t)\big)c^2\Big).
$$
Inserting this inequality in (\ref{maxesp}), using the fact that $\vp_c(x)=c\vp_1(\sqrt{c}x)$ and integrating in $y$
give the bound
$$
\E \Big( \vp_{C^{\vep}(t)}(x-X^{\vep}(t)) \Big) \le \frac{K}{(\det \Sigma)^{1/2}} \Int_0^{+\infty} \sqrt{c+c_0}
e^{-\frac12 \frac{\vep^2}{\det \Sigma} [ \s_{11} \frac{t^3}{12} +(\s_{22}-\frac{\s_{12}^2}{\s_{11}}t)]c^2} dc
$$
where $K$ is a constant, and since
$$
\Int_0^{+\infty} \sqrt{c} e^{-\frac{c^2}{2\a^2}} dc \le K \a^{3/2}
$$
for another constant $K$, it follows
\begin{equation}
\label{diffusion}
\max_{x\in \R} \E \Big( \vp_{C^{\vep}(t)}(x-X^{\vep}(t)) \Big) \le K_0 \vep^{-1/2} t^{-5/4}
\end{equation}
for $t$ large enough.

This inequality has to be compared to the result of \cite{Wa} where an additive equation 
with a white noise in time was considered. An inequality of the form (\ref{diffusion})
was obtained, but with a power $t^{-3/2}$ instead of $t^{-5/4}$.

\end{document}